\documentclass[preprint,3p,12pt]{elsarticle}

\usepackage{amssymb}

\usepackage{lipsum}
\usepackage{subfigure}
\usepackage{amsfonts}
\usepackage{amsthm}
\usepackage{nicefrac}
\usepackage{graphicx}
\usepackage{epstopdf}
\usepackage{amsfonts}
\usepackage{graphicx}
\usepackage{amsopn}
\usepackage{listings}
\usepackage{xcolor}
\usepackage{float}
\usepackage{marvosym}
\usepackage{booktabs}
\usepackage{multirow}
\usepackage{algorithm}
\usepackage{algorithmic}
\usepackage{graphicx}
\usepackage{amsmath}
\usepackage{amsfonts}
\usepackage{amsbsy}
\usepackage{amscd}
\usepackage{amssymb}
\usepackage{latexsym}
\usepackage[mathscr]{eucal}
\usepackage{bm}
\usepackage{threeparttable}
\usepackage{rotating}
\usepackage{tikz}

\usepackage{lineno}
\modulolinenumbers[1]

\definecolor{mygray}{gray}{0.9}

\newfont{\bb}{msbm10}

\def\vecx{\mathbf{x}}
\def\vecy{\mathbf{y}}

\def\vecr{\mathbf{r}}
\def\vecf{\mathbf{f}}
\def\vecres{\mathbf{res}}

\def\vecu{\mathbf{u}}
\def\vecv{\mathbf{v}}

\def\vecs{\mathbf{s}}

\def\matA{\mathbf{A}}
\def\matB{\mathbf{B}}

\def\matM{\mathbf{M}}

\def\matU{\mathbf{U}}
\def\matV{\mathbf{V}}

\def\matX{\mathbf{X}}
\def\matY{\mathbf{Y}}

\theoremstyle{proposition}
\newtheorem{proposition}{Proposition}[section]

\renewcommand{\Re}{\operatorname{Re}}
\renewcommand{\Im}{\operatorname{Im}}

\journal{arXiv}

\begin{document}

\graphicspath{{figures/}} 


\begin{frontmatter}



\title{ParaSLRF: A High Performance Rational Filter Method for Solving Large Scale Eigenvalue Problems }



\author[a,c]{Biyi Wang}
\ead{wangbiyi20@gscaep.ac.cn}
\author[b]{Karl Meerbergen}
\author[b]{Raf Vandebril}
\author[c,d]{Hengbin An}
\author[c,d]{Zeyao Mo}
            
\affiliation[a]{organization={Graduate School of China Academy of Engineering Physics},
            city={Beijing},
            postcode={100088}, 
            country={China}}

\affiliation[b]{organization={Department of Computer Science, KU Leuven, Celestijnenlaan 200A},
            city={Leuven},
            postcode={3001}, 
            country={Belgium}}

\affiliation[c]{organization={Institute of Applied Physics and Computational Mathematics},
            city={Beijing},
            postcode={100094}, 
            country={China}}

\affiliation[d]{organization={CAEP Software Center for High Performance Numerical Simulation},
            city={Beijing},
            postcode={100088}, 
            country={China}}

\begin{abstract}

In \emph{Wang et al., A Shifted Laplace Rational Filter for Large-Scale Eigenvalue Problems}, the SLRF method was proposed to compute all eigenvalues of a symmetric definite generalized eigenvalue problem lying in an interval on the real positive axis.
The current paper discusses a parallel implementation of this method, abbreviated as ParaSLRF. 
The parallelization consists of two levels:
(1) on the highest level, the application of the rational filter to the various vectors is partitioned among groups of processors; (2) within each group, every linear system is solved in parallel.

In ParaSLRF, the linear systems are solved by iterative methods instead of direct ones, in contrast to other rational filter methods, such as, PFEAST.
Because of the specific selection of poles in ParaSLRF, the computational cost of solving the associated linear systems for each pole, is almost the same.
This intrinsically leads to a better load balance between each group of resources, and reduces waiting times of processes.

We show numerical experiments from finite element models of mechanical vibrations, and show a detailed parallel performance analysis.
ParaSLRF shows the best parallel efficiency, compared to other rational filter methods based on quadrature rules for contour integration.
To further improve performance, the converged eigenpairs are locked, and a good initial guess of  iterative linear solver is proposed.
These enhancements of ParaSLRF show good two-level strong scalability and excellent load balance in our experiments.
\end{abstract}




\end{frontmatter}



\section{Introduction}
\label{sect:introd}
Consider the following symmetric definite generalized eigenvalue problem (GEP)
\begin{eqnarray}
\label{eq:GEP}
\matA \vecx = \lambda \matB \vecx, 
\end{eqnarray}
where $\matA \in \mathbb{R}^{n \times n}$ and $\matB \in \mathbb{R}^{n \times n}$
are large, sparse and symmetric positive definite matrices.
The GEP~\eqref{eq:GEP} arises from, e.g.,  vibration analysis, quantum mechanics, electronic structure calculations, etc.~\cite{saad2011numerical}.
For example, in vibration analysis, the matrix $\matA$ is the stiffness matrix and $\matB$ is 
the mass matrix.
The tuple $(\lambda, \vecx)$ is called an eigenpair of the GEP.
Because of the positive definiteness of the matrices $A$ and $B$, the eigenvalues are real and positive.
We are interested in computing all eigenvalues in a specified interval $(0, \gamma]$ and
their associated eigenvectors.

To achieve this goal, there is a class of methods based on contour integration, utilizing the resolvent operator along a closed contour, which encloses all the desired eigenvalues.
This class of methods was proposed in 2003 by Sakurai and Sugiura~\cite{sakurai2003projection,ikegami2010filter}.
The approximation of the contour integral by a quadrature rule is a \emph{rational filter function}.
There are several quadrature rules  that can be used, such as  midpoint~\cite{xi2016computing}, Gauss-Legendre~\cite{polizzi2009density}, Gauss-Chebyshev~\cite{xi2016computing}, and Zolotarev~\cite{guttel2015zolotarev} quadrature rules.
The resulted rational filter functions can amplify eigenvalues within the contour and suppress eigenvalues outside the contour.
In addition to the rational filter functions obtained by contour integrals, 
Xi and Saad~\cite{xi2016computing} and Wang, Meerbergen, Vandebril, An, and Mo \cite{wang2025SLRF}  determine the rational filter function differently; they utilize complex conjugate poles and put the poles on lines through the origin with a particular slope.

The application of the rational filter function requires the solution of $\ell$ systems with multiple right-hand sides of the form
\begin{equation}
    (\matA - \sigma_j \matB) \matY = \matB\matV 
    \quad \text{for} \quad j = 1, \ldots, \ell,
    \label{equ: linear systems for rational function}
\end{equation}
where $\sigma_j$, $j=1,\ldots,\ell$, are the poles.
The choice of the poles is not only important for the quality of the filter function, but also has significant impact on the performance of the linear system solver.
The computational cost of solving~\eqref{equ: linear systems for rational function} with a direct method is almost independent of $\sigma_j$.
When the problem is large, iterative methods may become preferable.
In this case, the difficulty and cost for solving \eqref{equ: linear systems for rational function} strongly depends on the value of $\sigma_j$.
This paper only considers the case where the linear systems are solved by iterative methods.

Matrices~\eqref{eq:GEP} originating from mechanical vibration problems have properties similar to matrices in the Helmholtz equation.
In this paper, instead of using the poles based on various quadrature rules, we select the poles with the form of the shifted Laplace preconditioner for the Helmholtz equation~\cite{erlangga2004class,erlangga2006novel,van2007spectral}, which was inspired by 
earlier work
\cite{magolu2000preconditioning,magolu2001incomplete,made2004performance}.
The main idea is to pick  complex conjugate poles 
on the lines $\{y = x(1 \pm \alpha \imath), \alpha > 0, x\in\mathbb{R}^+\}$.
We call the corresponding rational filter a shifted Laplace rational filter (SLRF)~\cite{wang2025SLRF}. It makes the linear system complex, but easier to precondition and solve. A detailed analysis and discussion is given by Wang et al.~\cite{wang2025SLRF}.

To the best of our knowledge, no existing study has thoroughly investigated the parallel performance of rational filter methods when the corresponding linear systems are solved using iterative linear solvers.
In this paper, we experimentally show and analyze the impact of different strategies on the performance of various parallel rational filter methods when computing many eigenpairs.
Our results show that there is a severe work-load imbalance and huge communication waiting overhead for the rational filters based on the classical quadrature rules.
The parallel shifted Laplace rational filter (ParaSLRF) performs best in terms of both CPU time and computation efficiency in our experiments.
To further reduce the computational cost, we use a locking strategy and reduce the number of vectors required to be filtered, when convergence is reached.
In addition, we propose to use a scaled Ritz vector as initial guess to accelerate the convergence of the linear system~\eqref{equ: linear systems for rational function}.
To reduce the computation cost further and achieve a better load-balancing, we fix the maximum number of linear iterations.

The outline of this paper is as follows. 
In Section~\ref{sect: Rational filter}, we review parallel rational filter methods for eigenvalue computations
and briefly discuss our strategy of selecting the poles and weights for the rational filter. 
We compare the computational cost between ParaSLRF and the rational filters based on various quadrature rules.
Additionally, we illustrate the data distribution patterns and the computation flow within the parallel implementation framework.
In Section~\ref{sect:numer-result}, we present the numerical results of the ParaSLRF and compare its computational efficiency with that of the aforementioned filters based on quadrature rules.
In Section~\ref{sect: enhanced ParaSLRF}, we adopt two strategies to improve the performance further and represent strong scalability of the enhanced version.
We end this paper with the main conclusions in Section~\ref{sect:conclusion}.

The following notation is used throughout this paper:
$\matA \in \mathbb{R}^{n \times n}$ and $\matB \in \mathbb{R}^{n \times n}$
represent $n \times n$ large, sparse and symmetric positive definite matrices.
The $\matB$-inner product of two vectors $\vecu$, $\vecv \in \mathbb{R}^{n}$ is defined by
$
\langle \vecu, \vecv \rangle_{\matB} =
\vecv^{\top} \matB \vecu.
$
The induced norm is $\| \vecv \|_{\matB} = \sqrt{\langle \vecv, \vecv \rangle_{\matB}}$.
We define $N$ as the total number of poles $\sigma$ (in the upper plane), with corresponding weights denoted by $w$.

\section{Parallel rational filter method}
\label{sect: Rational filter}
The core of the rational filter method is the rational function
\begin{equation}
\Phi (\lambda) = \sum_{j=1}^{\ell} \frac{w_j}{\lambda - \sigma_j},
\label{equ: rational function filter}
\end{equation}
where $\sigma_j$ is the \textit{j}-th pole, and $w_j$ is the corresponding weight, for $j=1,2,\ldots,\ell$.
For a closed contour $\Gamma\subset\mathbb{C}$ which encloses all the desired eigenvalues, the rational function~\eqref{equ: rational function filter} can amplify the function value at the desired eigenvalues $\lambda_i$ inside the contour $\Gamma$, and decrease the function values at the eigenvalues outside the contour $\Gamma$.
The name \textit{rational filter} arises from the filter properties.

To apply the filter it is necessary to solve $\ell$ linear systems of the form~\eqref{equ: linear systems for rational function}.
Since $\matA$ and $\matB$ are real symmetric and the domain of interest is an interval on the real axis, we choose the $\ell$ nodes as complex conjugate points in the complex plane. Let $\ell=2N$, then, with a proper ordering, we have $\sigma_{N+j}=\overline{\sigma_j}$ and $w_{N+j}=\overline{w_j}$, $j=1,\ldots,N$.
As a result, $N$ complex valued linear systems have to be solved when the filter is applied to a real vector.
The computation can be done in parallel.

The main idea is to distribute the required linear systems to multiple computing resources and solve them in parallel, then, use summation and reduction operations.
We refer to each rational filtering step as an outer iteration, and the linear solves within each filtering step as inner iterations.
For brevity, we assume the $np$ processes and $N$ poles  partitioned in $n_{part}$ groups.
We create parallel sub-communicators $\{ G_j \}^{n_{part}}_{j=1}$ with $N_{sub}=N/n_{part}$ assigned poles and $np_{sub}=np/n_{part}$ processes.
A framework of the parallel rational filter method for computing eigenpairs of the matrix pencil~\eqref{eq:GEP} is given in Algorithm~\ref{alg:Parallel Rational filter method}, with details below. The algorithm aims to compute $\mbox{NEV}$ eigenvalues in $(0,\gamma]$ and associated eigenvectors.

\begin{algorithm}[htbp]
	\caption{Parallel rational filter method for GEP~\eqref{eq:GEP}}
	\label{alg:Parallel Rational filter method}
	\begin{algorithmic}[1]
		\REQUIRE $\matA \in \mathbb{R}^{n\times n}$, $\matB \in \mathbb{R}^{n\times n}$,  interval $(0, \gamma]$, $\mbox{NEV}$ wanted eigenvalues, $L$ columns of the starting vectors $\matV_1 \in \mathbb{R}^{n\times L}$ where $L \geq \mbox{NEV}$, 
        poles $\sigma_1,\ldots,\sigma_N \in \mathbb{C}$ and corresponding weights $w_1,\ldots,w_N\in \mathbb{C}$, $np$ number of processes, $n_{part}$ number of partitions.
        \\
		\ENSURE converged eigenpairs $[\Theta, \matX]$.\\
            \STATE  /* Initialization phase. */       
            \STATE Load matrix pencil $\matA, \matB$ and generate random starting vectors $\matV_1$ in a global communicator $G_0$.
            \STATE Split the global communicator $G_0$ into $n_{part}$ sub-communicators $G_j, j=1,2,\ldots,n_{part}$, where the number of assigned processes and poles are $np_{sub}=np/n_{part}$ and $N_{sub}=N/n_{part}$, respectively. 
            Scatter and redistribute $\matA, \matB$ and $\matV_1$ from the global communicator to each sub-communicator.
		\FOR{$k = 1,2,\ldots$ until all of the wanted eigenpairs converged}
                \STATE  /* Parallel over the $n_{part}$ partitions. */
                \FOR{$j = 1,2,\ldots$, $n_{part}$ \textbf{in parallel}}
                    \FOR{$i = 1,2,\ldots$, $N_{sub}$}
                        \STATE Set $\sigma_{j,i} = \sigma_{(j-1)N_{sub}+i}$ and $w_{j,i} = w_{(j-1)N_{sub}+i}$.
                        \STATE Solve $\matY_{j,i} = (\matA - \sigma_{j,i} \matB)^{-1}\matB\matV_k$ in sub-communicator $G_j$ \textbf{in parallel}.
                    \ENDFOR
                \ENDFOR
    \STATE Concatenate $\matY_j=[\matY_{j,1},\matY_{j,2},\ldots,\matY_{j,N_{sub}}]$ in sub-communicator $G_j$, \textbf{for} $j=1,\ldots,n_{part}$ \textbf{in parallel}. 
    \STATE Scatter $\matY_j$ from each sub-communicator to the global communicator and then perform a sum-reduce operation $\matU = \sum_{j=1}^{n_{part}} \sum_{i=1}^{N_{sub}} 2\Re(w_{j,i} \matY_{j,i})$.
                \STATE Solve the projected eigenvalue problem $(\matU^\top \matA \matU) \vecs_m = \theta_m (\matU^\top \matB \matU) \vecs_m$ and let $\vecv_m=\matU \vecs_m$, $m=1,2,\ldots,L$, with $\theta_m$ sorted in ascending order. 
                \STATE Check the convergence of the approximated eigenpairs $(\theta_m,\vecv_m),m=1,2,\ldots,L$.
                If converged, append $\theta_m$ to $\Theta$, and duplicate $\vecv_m$ to $\matX$, respectively.
                \STATE Let $\matV_{k+1}=[\vecv_1,\vecv_2,\ldots,\vecv_L]$.
                Scatter $\matV_{k+1}$ to each sub-communicator again.
            \ENDFOR
	\end{algorithmic}
\end{algorithm}

On Lines~2 to 3, the initialization is executed. 
The matrices $\matA, \matB, \matV_1$ are duplicated and redistributed across each sub-communicator $G_j$.
On Lines~6 to 11, the computing tasks are divided into $n_{part}$ parts and the computation is performed in parallel, in multiple sub-communicators.
On Line~13, the matrix $\matU$ calculated by the sum-reduce operation represents the vectors obtained by applying the rational filtering on vectors $\matV_k$.
On Line~14, the projected pencil consisting of matrices $\hat{\matA} = \matU^{\top} \matA \matU$ and $\hat{\matB} = \matU^{\top} \matB \matU$ is formed.
The approximated eigenpairs $(\theta_m,\vecv_m)$ are extracted from the subspace projection on $\matA$ and $\matB$, by using the Rayleigh-Ritz method. 
On Line~15, the convergence of the approximated pairs $(\theta_m,\vecv_m)$ is examined.
On Line~16, the matrix $\matV_{k+1}$ is assembled and then scattered to each sub-communicator for the next rational filtering iteration.
As described in Algorithm~\ref{alg:Parallel Rational filter method}, we can use a two-level parallelism: 
\begin{itemize}
    \item Level-1: solve the linear systems with multiple right-hand terms independently, over $n_{part}$ sub-communicators $G_j$, $j=1,\ldots,n_{part}$;
    
    \item Level-2: solve each linear system in parallel within sub-communicator $G_j$.
\end{itemize}

\subsection{Data distribution and computation flow}
\label{subsec: data structure}

\begin{figure}[htbp]
\centering
\includegraphics[scale=0.25]{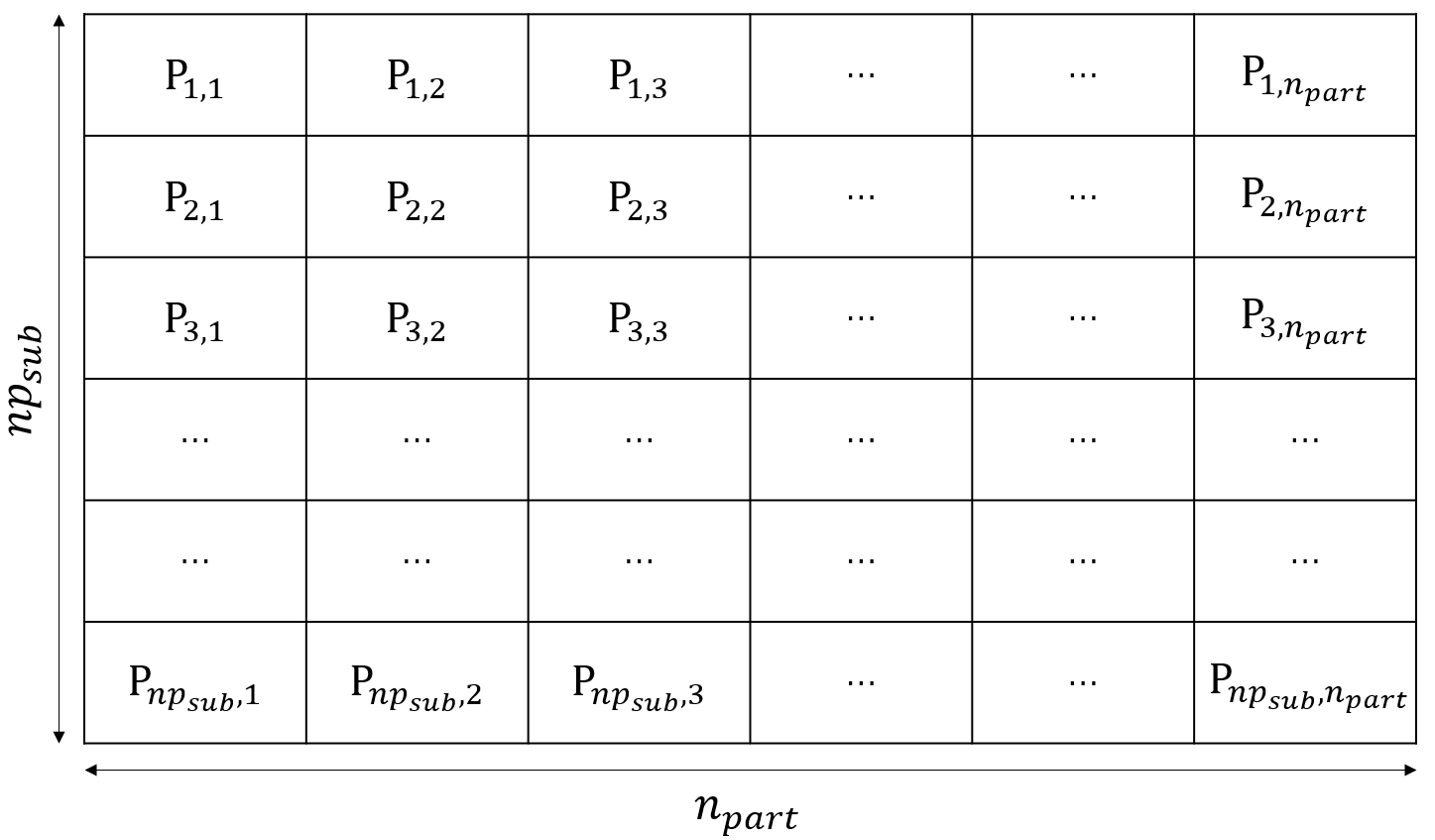}
\caption{The 2D representation of $np$ processes.}
\label{fig: 2D grid}
\end{figure}

To achieve two-level parallelism , we need to manage the data properly and design a suitable computation flow.
Assume we have $np$ processes, $N$ poles in total. We partition the $np$ processes into $n_{part}$ groups. (For the ease of presentation, we assume here that $np$ is a multiple of $n_{part}$.)
By default, $np$ processes are organized as a 1D vector to represent the global communicator \texttt{MPI\_COMM\_WORLD}, denoted by $G_{0}$.
The $np$ processes are reorganized into a 2D grid of size $n_{part} \times (np / n_{part})$.
Each column $P(:,j)$, where $ j=1,\ldots,n_{part}$, represents a group of resources. We associate  it with a sub-communicator $G_j$ consisting of $np_{sub}=(np / n_{part})$ processes, see Figure~\ref{fig: 2D grid}.
For each column $j$, we assign $N_{sub}=N/n_{part}$ poles.

\begin{figure}[htbp]
\centering
\includegraphics[scale=0.45]{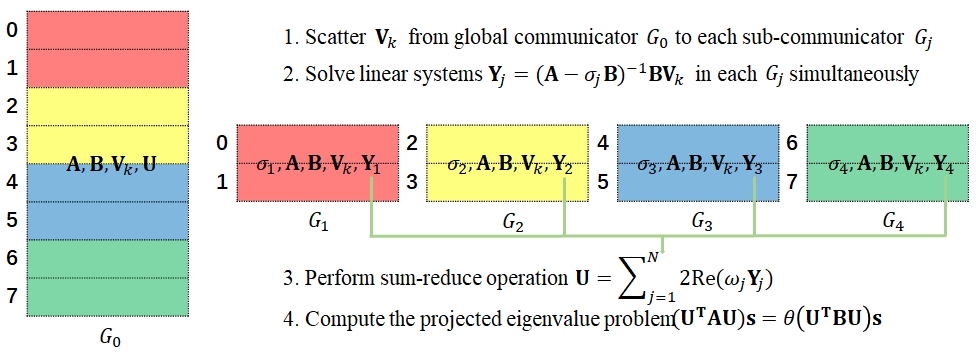}
\caption{Computational flow of the parallel rational filter method at iteration step $k$, for $np=8, N=n_{part}=4$.}
\label{fig:computationflow}
\end{figure}

First we load the matrix pencil to the  global communicator $G_{0}$, and each process stores some consecutive matrix rows.
In $G_{0}$, we generate $L$ random  vectors $\matV_1$ with the same row distribution as the matrix pencil.
Then we scatter the matrix pencil and  vectors forward to each sub-communicator $G_{j}$. 
We follow the same data storage and manipulation as in SLEPc~\cite{techreport_CISS}.
To illustrate Algorithm~\ref{alg:Parallel Rational filter method} more clearly,  Figure~\ref{fig:computationflow}  shows the data distribution and detailed operations at iteration step  $k$.

In this example $np=8$ and $N=4$, and the global communicator $G_0$ is split into 4 parts, i.e., each sub-communicator $G_j$ manages 2 processes.
For each $G_j$, we store the full $\matA, \matB$, and assign one pole $\sigma_j$.
The right-hand-side vectors $\matV_k$ maintain the same row-split distribution as the matrix pencil, consistently preserved across both the global communicator and all sub-communicators.
After solving the required linear systems in each $G_j$, the solution vectors $\matY_j$ are scattered backward to the global communicator, where the sum-reduce operation is performed to compute the rational filtered vectors $\matU$.
Then, we form the projected pencil $(\matU^{\top}\matA\matU,\matU^{\top}\matB\matU)$, where $\matU, \matA$ and $\matB$ have the same row-split distribution in the global communicator.
As the projected problem is small, we store the matrix product in the root process of $G_0$.
Finally, we solve the eigenpairs of the small projected problem, which are used for updating the approximated vectors $\matV_{k+1}$ for next iteration step.

As we can see, there is no explicit orthogonalization operation, so there is no communication amongst the $G_j$'s, i.e., groups of processes corresponding to different poles.
However, the sum-reduce operation requires all the linear systems to be solved in each sub-communicator. 
When the linear systems assigned in each sub-communicator have different difficulty levels, there might be waiting time in the MPI communication between the sub-communicators.

\subsection{The choice of poles}
\label{subsec: optimal poles}
The main difference between the various rational filter methods is the selection of poles $\sigma_j$, and corresponding weights $w_j$.
Apart from the classical quadrature rules, there are other interesting ways to select the poles and weights~\cite{xi2016computing,guttel2015zolotarev,austin2015computing,ohno2010quadrature}.

To solve the associated linear systems for each pole, iterative solvers are preferable for large-scale problems.
However, we have to take the position of the poles into account carefully.
Given a pole $\sigma$ we can represent it as 
\[
\sigma=\Re(\sigma)(1 + \alpha \imath),
\]
where $\alpha$ equals $\Re(\sigma)/\Im(\sigma)$ and indicates the degree of deviation from the real axis.
The effect of $\alpha$ on solving the associated linear systems is analysed by Wang, et al.~\cite{wang2025SLRF}.
On the one hand, for example, part of the poles generated from Gauss-Chebyshev quadrature, are very close to the real axis (i.e., small $\alpha$) usually result in ill-conditioned linear systems, which are more difficult to solve by iterative solvers.
On the other hand, as pointed out in~\cite{meerbergen1997restarted,ohno2010quadrature,xi2016computing,wang2025SLRF}, moving poles away from the spectrum (or real axis for Hermitian problem) can reduce the difficulty of solving the associated linear systems by using an iterative solver.
In this paper, the poles and weights from the shifted Laplace rational filter (SLRF) are determined with parameters $\alpha = 1$ and $\beta = 0.01$ which is a relaxation parameter~\cite{wang2025SLRF}.

\section{Numerical experiments}
\label{sect:numer-result}
\subsection{Finite element model of vibration}
\label{subsect: eigen-problem}
We consider the following vibration equation
\begin{eqnarray}
\label{eq:governing equation}
\left\{
\begin{array}{rcl}
\bm{\sigma} \cdot \nabla -\rho \ddot{\mathbf{u}} &=& 0, \quad \text{in} \; \Omega, \\
\bm{\sigma} \cdot \mathbf{n} = \mathbf{\overline{t}}  &=& 0, \quad \text{on} \; \Gamma_N, \\
\mathbf{u}                                       &=& 0, \quad \text{on} \; \Gamma_D,
\end{array}
\right.
\end{eqnarray}
where $\mathbf{u}$ is the displacement, 
$\bm{\sigma}$ is the stress and 
$\rho$ is the density of the material.
The domain is denoted by $\Omega$, 
$\Gamma_D$ denotes the Dirichlet boundary condition and 
$\Gamma_N$ denotes the Neumann boundary condition, 
$\mathbf{n}$ is the outer normal vector of the boundary of the domain $\partial \Omega$.
The derivation of the weak formulation of~\eqref{eq:governing equation} for a finite element discretization is discussed by Wang, An, Xie, and Mo~\cite{wang2025new}.

In the following, we introduce three mechanical models for numerical tests.
These models have the same material parameters: 
the density of the material is $1.0\text{kg}/\text{m}^2$,
the Poisson's ratio is 0.29,
the Young's modulus is $21.5\text{N}/\text{m}^2$.

\subsubsection{Pyramid model}
\label{subsect: pyramid model}
\begin{figure}[htbp]
\centering
\includegraphics[scale=0.25]{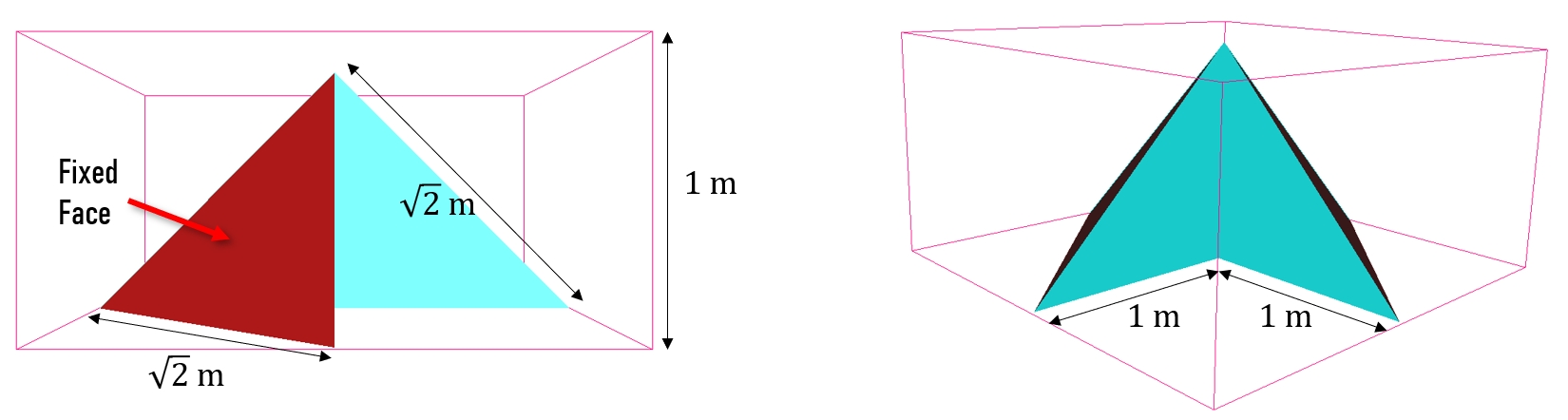}
\caption{The shape of pyramid model.}
\label{model:pyramid}
\end{figure}

The first model is a three dimensional mechanical model which looks like a pyramid.
The foundation is fixed on the ground and two surfaces in red are fixed, see Figure~\ref{model:pyramid}. 
The domain is partitioned by tetrahedral elements and the equations are discretized by first order
finite elements.
Two scales of this problem will be tested:
\begin{itemize}
    \item Pyramid S1: $n = 36,339$ with $12,113$ grid nodes; 
    \item Pyramid S2: $n = 268,353$ with $89,451$ grid nodes; 
\end{itemize}
where $n$ is the number of the DOFs.

For Pyramid S1, the first 20 smallest eigenvalues are located in [350.6497, 1530.2328], and the first 100 smallest eigenvalues are located in [350.6497, 4272.0205].
For Pyramid S2, the first 20 smallest eigenvalues are located in [344.4414,1477.9992], and the first 100 smallest eigenvalues are located in [344.4414,3907.2456].

\subsubsection{Hollow platform model and fish-like model}
\label{subsect: hollow platform model and fish-like model}

The hollow platform model originates from cutting out a cylinder 
from a platform. 
The surface of the cylinder is fixed. 
The dimension of the obtained GEP is 37161.
The fish-like model is generated by rotating a given two dimensional region along the $x$-axis. 
The obtained body looks like a fish and the lower surface is fixed.
The dimension of the obtained GEP is 55248.
For more details about this two models, we refer to Wang, et al.
 \S 4 in~\cite{wang2025SLRF}.

\subsection{Setting}
\label{subsect: setting}
The computations were carried out on a Linux compute server
with two 28-cores Intel Xeon Gold 6330N processors, and 1 TB of RAM.
All methods are written in C, 
built on top of PETSc~\cite{petsc} and SLEPc~\cite{slepc}, and compiled with the \texttt{-O3} optimization flag.
For a fair comparison, the same matrix-vector products, inner products interfaces integrated in PETSc and SLEPc are called for all methods. 
In the following, we define $\vecr_k = \frac{(\matA - \theta_k \matB) \vecx_k}{\theta_k \Vert\vecx_k \Vert_{\matB}}$
as the relative residual vector where $\theta_{k}$
is Rayleigh quotient associated with the approximate eigenvector $\vecx_{k}$.
As stopping criteria, we use $\| \vecr_{k}\|_2 < 10^{-8}$.

We have to point out that the incomplete LU factorization method (ILU) and BoomerAMG method integrated in HYPRE~\cite{falgout2002hypre} lack robust support for complex linear system computations.
Unless otherwise stated, the linear solver is GCR (generalized conjugate residual) method with the GAMG preconditioner in PETSc.
We reuse the preconditioners for each pole to avoid assembling the associated matrices repeatedly during the construction phase of GAMG.
For the linear system $\matM \vecu =\vecf$, the stopping criteria for the linear solver requests the relative residual norm  $\| \vecf - \matM \vecu_k\|_2/ \| \vecf\|_2$ to be less than $10^{-10}$ or the number of iterations has reached the maximum iteration number, set to 10000.

For the parallel shifted Laplace rational filter (ParaSLRF) the values of poles and corresponding weights are determined by the parameters $\alpha = 1$ and $\beta = 0.01$. 
The dimension of the search space $L$ should be larger than the number of desired eigenvalues ($\mbox{NEV}$). Therefore, we choose $L=1.2\times\mbox{NEV}$.

\subsection{Results and efficiency analysis}
\label{subsect: results}
The first 20 and 100 smallest eigenpairs of the aforementioned mechanical models are computed.
For each computation task we use $np=16$ processes.
For each rational filter method we use $N=4$ poles (in the upper plane).
The processes are split into $n_{part}=4$ parts, i.e., each sub-communicator is bound to one pole and manages $np_{sub}=np/n_{part}=4$ processes for solving the associated linear systems. 

In Table~\ref{tab: numerical results}, the reported results include the number of outer iterations \textbf{Iter} and the CPU time \textbf{Time$_{\text{Total}}$}. Moreover, we record the maximum time, denoted by \textbf{Time$_{\text{LinearSolve}}^{\text{Max}}$} for solving the required linear systems across all sub-communicators and we document its ratio, denoted by \textbf{Ratio}, relative to the minimum time of linear solves in the whole procedure.

\begin{table}[htbp]
  \centering
  \caption{Numerical results.}
  \begin{threeparttable}
    \begin{tabular}{c|cccccc}
    \multicolumn{1}{c}{Model} & Interval & Filter & Iter  & Time$_{\text{Total}}$  & Time$_{\text{LinearSolve}}^{\text{Max}}$ & Ratio \\
    \midrule
    \multicolumn{1}{c|}{\multirow{8}[4]{*}{Pyramid S1}} & \multicolumn{1}{c}{\multirow{4}[2]{*}{\shortstack{NEV=20 \\ (0, 1531]}}} & Midpoint &   15    &   2.0292E+03 & 2.0279E+03    &  4.9 \\
          &       & Gauss-Legendre &   5    &   1.2217E+03 & 1.2210E+03     &  9.0 \\
          &       & Gauss-Chebyshev &   9    &   4.0055E+03 & 4.0045E+03     & 16.3 \\
          &       & ParaSLRF &   10    &   3.2080E+02 & 3.1960E+02   & 1.0 \\
\cmidrule{2-7}          & \multicolumn{1}{c}{\multirow{4}[2]{*}{\shortstack{NEV=100 \\(0, 4273]}}} & Midpoint &  10    &  8.2540E+03 & 8.2484E+03   & 7.2 \\
          &       & Gauss-Legendre &   5    &   7.8351E+03 & 7.8322E+03  & 13.5 \\
          &       & Gauss-Chebyshev &   9    &  2.7157E+04 & 2.7152E+04  & 26.0 \\
        &       & ParaSLRF &   9    &  1.1535E+03 & 1.1484E+03   & 1.0 \\
    \midrule
    \multirow{8}[4]{*}{\shortstack{Hollow \\ platform}} & \multicolumn{1}{c}{\multirow{4}[2]{*}{\shortstack{NEV=20 \\ (0, 215]}}} & Midpoint &  4     &  4.2453E+03 & 4.2434E+03    & 4.7 \\
          &       & Gauss-Legendre &   4    &  7.7003E+03 & 7.6985E+03     & 8.6 \\
          &       & Gauss-Chebyshev &  4     &  1.4001E+04 & 1.3999E+04     &  15.9 \\
          &       & ParaSLRF &  6     &   1.5403E+03 & 1.5382E+03     & 1.0 \\
\cmidrule{2-7}          & \multicolumn{1}{c}{\multirow{4}[2]{*}{\shortstack{NEV=100 \\ (0, 1507]}}} & Midpoint &  7    &  3.5463E+04 & 3.5456E+04      & 6.4 \\
          &       & Gauss-Legendre &  4     &  4.0816E+04 & 4.0812E+04     & 12.8 \\
          &       & Gauss-Chebyshev &  5     &  1.0494E+05 & 1.0493E+05    & 25.8 \\
          &       & ParaSLRF &   7    &    5.4263E+03 & 5.4121E+03     & 1.0 \\
    \midrule
    \multirow{8}[4]{*}{Fish-like} & \multicolumn{1}{c}{\multirow{4}[2]{*}{\shortstack{NEV=20 \\ (0, 1858]}}} & Midpoint &  6     &  6.3236E+03 & 6.3226E+03     & 4.1 \\
          &       & Gauss-Legendre &   5    &  9.8096E+03 & 9.8088E+03      & 7.8 \\
          &       & Gauss-Chebyshev &  6     &  2.1733E+04 & 2.1733E+04     &  14.3 \\
          &       & ParaSLRF &  7     &   1.6070E+03 & 1.6061E+03    & 1.1 \\
\cmidrule{2-7}          & \multicolumn{1}{c}{\multirow{4}[2]{*}{\shortstack{NEV=100 \\ (0, 10239]}}} & Midpoint &   11    &   3.9331E+04 & 3.9324E+04      & 6.2 \\
          &       & Gauss-Legendre &  6     &  4.0080E+04 & 4.0076E+04    & 12.5 \\
          &       & Gauss-Chebyshev &  9     &  1.1449E+05 & 1.1449E+05    & 24.2 \\
          &       & ParaSLRF &   7    &    3.8386E+03 & 3.8346E+03    & 1.0 \\
    \bottomrule
    \end{tabular}%
    \end{threeparttable}
  \label{tab: numerical results}%
\end{table}%

From Table~\ref{tab: numerical results}, we observe that ParaSLRF achieves the best performance in our experiments, it successfully completes all tasks with the shortest execution time.
Instead, the rational filters based on the quadrature rules do not work well. 
For instance, the Gauss-Chebyshev filter demonstrates the highest computational cost, requiring up to 29.8 times more CPU time compared to the ParaSLRF method, for solving the first 100 eigenpairs of the fish-like model.
From the penultimate column, we see that solving the associated linear systems is the most expensive part in the whole procedure.
However, we have to point out that the allocated processes are not fully utilized and many processes are idle most of the time, for the rational filters based on quadrature rules.
From the last column of Table~\ref{tab: numerical results}, we can observe that the \textbf{Ratio} value can be really high which means severe imbalance between the sub-communicators.
The reason for this phenomenon is that the time for solving the linear systems in each sub-communicator is quite different.
The reduction of vectors $\matU$ requires  synchronization,
forcing faster sub-communicators to wait for the slowest one to complete the linear solve.
On the contrary, ParaSLRF maintains a consistent \textbf{Ratio} value approximately equal to 1, indicating that 
(i) the computational cost for solving the associated linear systems is almost constant across all sub-communicators;
(ii) the utilization of resources is high and the MPI waiting overhead is low.

To validate our statements, 
we record the additional MPI waiting time (i.e.\ the maximum idle time across all sub-communicators), denoted by \textbf{Time$_{\text{MPIWait}}^{\text{Max}}$}, during the stage of the sum-reduce operation and report its proportion to the total CPU time, denoted by \textbf{Prop}, for each task. The results are represented in Table~\ref{tab: MPI waiting}.

\begin{table}[htbp]
  \centering
  \caption{Extra MPI waiting time and its proportion to total solution time.}
  \begin{threeparttable}
    \begin{tabular}{c|cccc}
    \multicolumn{1}{c}{Model} & Interval & Filter & Time$_{\text{MPIWait}}^{\text{Max}}$ & Prop (\%) \\
    \midrule
    \multicolumn{1}{c|}{\multirow{8}[4]{*}{Pyramid S1}} & \multicolumn{1}{c}{\multirow{4}[2]{*}{\shortstack{NEV=20 \\ (0, 1531]}}} & Midpoint &  1.6182E+03   & 79.75 \\
          &       & Gauss-Legendre &   1.0849E+03    & 88.80 \\
          &       & Gauss-Chebyshev &  3.7596E+03    & 93.86 \\
          &       & ParaSLRF &  8.1755E+00    & 2.52 \\
\cmidrule{2-5}          & \multicolumn{1}{c}{\multirow{4}[2]{*}{\shortstack{NEV=100 \\(0, 4273]}}} & Midpoint &   7.1097E+03  & 86.14 \\
          &       & Gauss-Legendre &   7.2531E+03    & 92.57 \\
          &       & Gauss-Chebyshev &  2.6109E+04      & 96.14 \\
          &       & ParaSLRF &  1.0824E+01   & 0.94 \\
    \midrule
    \multirow{8}[4]{*}{\shortstack{Hollow \\ platform}} & \multicolumn{1}{c}{\multirow{4}[2]{*}{\shortstack{NEV=20 \\ (0, 215]}}} & Midpoint & 3.3365E+03   & 78.59 \\
          &       & Gauss-Legendre &  6.8059E+03   & 88.38 \\
          &       & Gauss-Chebyshev & 1.3119E+04    & 93.70 \\
          &       & ParaSLRF &  3.4959E+01  & 2.27 \\
\cmidrule{2-5}          & \multicolumn{1}{c}{\multirow{4}[2]{*}{\shortstack{NEV=100 \\ (0, 1507]}}} & Midpoint &  2.9926E+04   & 84.39 \\
          &       & Gauss-Legendre &   3.7624E+04    & 92.18 \\
          &       & Gauss-Chebyshev &  1.0086E+05    & 96.11 \\
          &       & ParaSLRF & 1.2405E+02   & 2.29 \\
    \midrule
    \multirow{8}[4]{*}{Fish-like} & \multicolumn{1}{c}{\multirow{4}[2]{*}{\shortstack{NEV=20 \\ (0, 1858]}}} & Midpoint &  4.7763E+03      & 75.53 \\
          &       & Gauss-Legendre & 8.5520E+03  & 87.18 \\
          &       & Gauss-Chebyshev & 2.0209E+04   & 92.99 \\
          &       & ParaSLRF  &  8.8895E+01    & 5.53 \\
\cmidrule{2-5}          & \multicolumn{1}{c}{\multirow{4}[2]{*}{\shortstack{NEV=100 \\ (0, 10239]}}} & Midpoint &  3.3013E+04    & 83.94 \\
          &       & Gauss-Legendre & 3.6881E+04  & 92.02 \\
          &       & Gauss-Chebyshev &  1.0976E+05   & 95.87 \\
          &       & ParaSLRF &  1.4755E+02   & 3.84 \\
    \bottomrule
    \end{tabular}%
    \end{threeparttable}
  \label{tab: MPI waiting}%
\end{table}%

From Table~\ref{tab: MPI waiting}, we can observe that the rational filter methods based on quadrature rules have high MPI waiting time.
The most evident example is the Gauss-Chebyshev rational filter with a maximum MPI waiting time of more than $92\%$ of the total time, that is, at least a quarter of the processes are idle more than $92\%$ of the time.
ParaSLRF on the other hand has a MPI waiting time of around $4\%$ of the total CPU time.
All computing resources are fully utilized as much as possible, because of the good work-load balancing.

\section{Enhanced ParaSLRF}
\label{sect: enhanced ParaSLRF}
As shown in Table~\ref{tab: numerical results}, solving the linear systems dominates the computational cost of the algorithm.
To further reduce the computational cost and improve the performance of ParaSLRF we adopt two useful strategies.

The first improvement consists of locking the approximate eigenpairs that satisfy the stopping criteria and reduce the number of vectors required to be filtered.
Next, in the subsequent iterations, we make the filtered vectors and the converged eigenvectors orthogonal with respect to the mass matrix $\matB$.

The second improvement lies in the fact that $\nicefrac{1}{(\theta-\sigma)} \vecv$, as the initial guess is a good choice to accelerate the linear solve, under the assumption that $(\theta,\vecv)$ is a good approximation of one of the eigenpairs $(\lambda,\vecx)$.
There are two arguments to see this.
First, employing an initial approximation close to the exact solution $\nicefrac{1}{(\lambda-\sigma)} \vecx$, seems an obvious strategy.
Second, it is known that a fixed number of iterations for the linear system with the Cayley transformation usually leads to an eigenvector approximation with a lower residual norm~\cite{meerbergen1997restarted,karl_Cayley,bai2000templates}; using the proposed initial guess for shift-and-invert transformation is actually applying the Cayley transformation with initial zero guess, as the following proposition shows.

\begin{proposition}
Given the pair $(\theta,\vecv)\in\mathbb{R}\times(\mathbb{R}^n\setminus\{0\})$. Let $\mathcal{K}^{T^C}_m(\matA - \sigma \matB,\vecres_C^{(0)})$ be the (preconditioned) Krylov subspace formed for solving the "Cayley transformation" linear system 
\begin{equation}
    T^C: 
    (\matA - \sigma \matB) \vecu = (\matA -\theta \matB) \vecv \quad \text{with  initial guess} \  \vecu^{(0)}=0,
    \label{equ: Cayley}
\end{equation}
and let $\mathcal{K}^{T^{SI}}_m(\matA - \sigma \matB,\vecres_{SI}^{(0)})$ be the Krylov space of dimension $m$ for solving the shift-and-invert transformation  system 
 \begin{equation}
    T^{SI}:
      (\matA - \sigma \matB) \vecy = \matB \vecv \quad \text{with  initial guess} \ \vecy^{(0)}=\frac{1}{\theta-\sigma} \vecv,
      \label{equ: shift-invert}
 \end{equation}
 then
\begin{itemize}
    \item[(i.)] the initial residuals satisfy  $\vecres^{(0)}_{SI} = \frac{1}{\sigma - \theta} \vecres^{(0)}_{C},$
    \item[(ii.)] the resulting Krylov subspaces $\mathcal{K}^{T^C}_m(\matA - \sigma \matB,\vecres_C^{(0)})$
    and $\mathcal{K}^{T^{SI}}_m(\matA - \sigma \matB,\vecres_{SI}^{(0)})$
    are identical.
\end{itemize}

\end{proposition}
\proof Split the right-hand side of~\eqref{equ: Cayley},
\begin{align*}
    (\matA - \sigma \matB) \vecu &= (\matA -\theta \matB) \vecv =(\matA - \sigma \matB) \vecv + (\sigma - \theta) \matB \vecv.
\end{align*}
thus, 
\begin{align*}
    (\matA - \sigma \matB) (\vecu - \vecv) &= (\sigma - \theta) \matB \vecv, \\
    (\matA - \sigma \matB)\frac{\vecu - \vecv}{\sigma - \theta} &= \matB \vecv.
\end{align*}
Let $\vecy=\frac{\vecu - \vecv}{\sigma - \theta}$, we know that $\vecu^{(0)}=0$ is equivalent to $\vecy^{(0)}=\frac{1}{\theta - \sigma}\vecv$. 
Let $P^{-1}$ denote the preconditioner.
For system~\eqref{equ: Cayley} the initial (preconditioned) residual is
\begin{equation*}
    \vecres^{(0)}_{C} = P^{-1}(\matA - \theta \matB) \vecv,
\end{equation*}
and for system~\eqref{equ: shift-invert} the initial residual is
\begin{equation*}
    \vecres^{(0)}_{SI} = P^{-1}\matB \vecv - P^{-1}(\matA - \sigma \matB) \frac{\vecv}{\theta - \sigma}=\frac{1}{\sigma - \theta}P^{-1}(\matA - \theta \matB) \vecv = \frac{1}{\sigma - \theta} \vecres^{(0)}_{C}.
\end{equation*}
The residuals are the starting vectors of the two Krylov spaces.
Since these vectors lie in the same direction and the (preconditioned) matrix is the same for both linear system, the two Krylov spaces must be the same.
\qed

The connection with the Cayley transform, shows that the scaled approximate eigenvector $\vecy^{(0)}=\frac{1}{\theta - \sigma}\vecv$ is a good initial guess of the iterative solver for the shift-and-invert linear system.
One can fix the maximum number of linear iterations $\rm It_{\text{max-linear}}$ of the linear solver, which is beneficial to reduce the cost of the linear solve and achieve a better load balance across all  sub-communicators.

The enhanced ParaSLRF is presented in Algorithm~\ref{alg: enhanced ParaSLRF}.
\begin{algorithm}[htbp]
	\caption{Enhanced ParaSLRF}
	\label{alg: enhanced ParaSLRF}
	\begin{algorithmic}[1]
		\REQUIRE $\matA \in \mathbb{R}^{n\times n}$, $\matB \in \mathbb{R}^{n\times n}$,  interval $(0, \gamma]$, 
        $\mbox{NEV}$ the number of the wanted eigenvalues, $L$ number of columns of the starting vectors where $L \geq \mbox{NEV}$,
        poles $\sigma_1,\ldots,\sigma_N \in \mathbb{C}$ and corresponding weights $w_1,\ldots,w_N\in \mathbb{C}$, $np$ number of processes, $n_{part}$ number of partitions, maximum linear iteration $\rm It_{\text{max-linear}}$.
        \\
		\ENSURE converged eigenpairs $[\Theta, \matX]$.\\
            \STATE  /* Initialization phase. */       
            \STATE Load matrix pencil $\matA, \matB$ and generate random starting vectors $\matV_1 \in \mathbb{R}^{n\times L}$ in global communicator $G_0$.
            \STATE Split the global communicator into $n_{part}$ sub-communicators where the number of assigned  processes are $np_{sub}=np/n_{part}$ and $N_{sub}=N/n_{part}$. 
            Scatter and redistribute $\matA, \matB$ and $\matV_1$ from the global communicator to each sub-communicator $G_j$.
            \STATE Set $n_{act}=L$ ($n_{act}$ counts the dimension of the active subspace), $n_{conv}=0$ ($n_{conv}$ counts the number of converged eigenpairs).
		\FOR{$k = 1,2,\ldots$ until all of the wanted eigenpairs have converged}
                \STATE  /* Parallel over the $n_{part}$ partitions. */
                \FOR{$j = 1,2,\ldots$, $n_{part}$}
                    \FOR{$i = 1,2,\ldots$, $N_{sub}$}
                        \STATE Set $\sigma_{j,i} = \sigma_{(j-1)N_{sub}+i}$ and $w_{j,i} = w_{(j-1)N_{sub}+i}$.
                        \STATE Solve $(\matA - \sigma_{j,i} \matB)\matY_{j,i} = \matB\matV_k$ by using the distributed linear iterative solver with  fixed maximum linear iterations $\rm It_{\text{max-linear}}$ in each sub-communicator $G_j$ at the same time.
                    \ENDFOR
                \ENDFOR
                \STATE Concat $\matY_j=[\matY_{j,1},\matY_{j,2},\ldots,\matY_{j,N_{sub}}]$ in each sub-communicator.
                \STATE Scatter each $\matY_j$ from sub-communicators to the global communicator where  a sum-reduce operation $\matU = \sum_{j=1}^{n_{part}} \sum_{i=1}^{N_{sub}} 2\Re(w_{j,i} \matY_{j,i})$ is performed. 
                \IF{$n_{conv}>0$}
                    \STATE Make the filtered vectors $\matU$ and converged eigenvectors $\matX$ orthogonal with respect to the mass matrix $\matB$.
                \ENDIF
                \STATE Solve the projected eigenvalue problem $(\matU^\top \matA \matU) \vecs_m = \theta_m (\matU^\top \matB \matU) \vecs_m$ and let $\vecv_m=\matU \vecs_m$, $m=1,2,\ldots,n_{act}$.
                \STATE Examine the convergence of each pair $(\theta_m,\vecv_m)$. If converged, append $\theta_m$ to $\Theta$ and duplicate $\vecv_m$ to $\matX$. Update the numbers $n_{conv}$, and $n_{act}=L-n_{conv}$.
                \STATE Remove the converged vectors in $\matV_{k}$ and resize $\matV_{k+1}=[\vecv_1,\ldots,\vecv_{n_{act}}]$.
                Scatter $\matV_{k+1}$ to each sub-communicator again.
                \STATE Let $\matY_{j,i}^{init}=[\frac{1}{\theta_1 - \sigma_{j,i}}\vecv_1,\ldots,\frac{1}{\theta_{n_{act}} - \sigma_{j,i}}\vecv_{n_{act}}], j=1,2,\ldots, n_{part}$, and $i=1,2,\ldots,N_{sub}$, be the initial guess of the iterative linear solver for the next iteration.
            \ENDFOR
	\end{algorithmic}
\end{algorithm}
On Line~10, the associated linear systems are solved by an iterative solver with the fixed maximum iteration number $\rm It_{\text{max-linear}}$.
On Lines~15 to 17, the filtered vectors $\matU$ which contain the information of the desired eigenvectors, are orthogonalized with the converged eigenvectors $\matX$, in order to avoid the converged eigenpairs from reappearing in the subsequent iterations.
On Lines~19 to 20, the converged eigenvectors are removed from $\matV_{k+1}$ so that only $n_{act}$ columns remain, where $n_{act}$ is the dimension of the active subspace for the next rational filtering iteration.
On Line~21, the initial guess $\matY_{j,i}^{init}$ is assembled by the scaled approximate eigenvectors to accelerate the next linear solve.

\begin{figure}[htbp]
	\centering
	\includegraphics[scale=0.4]{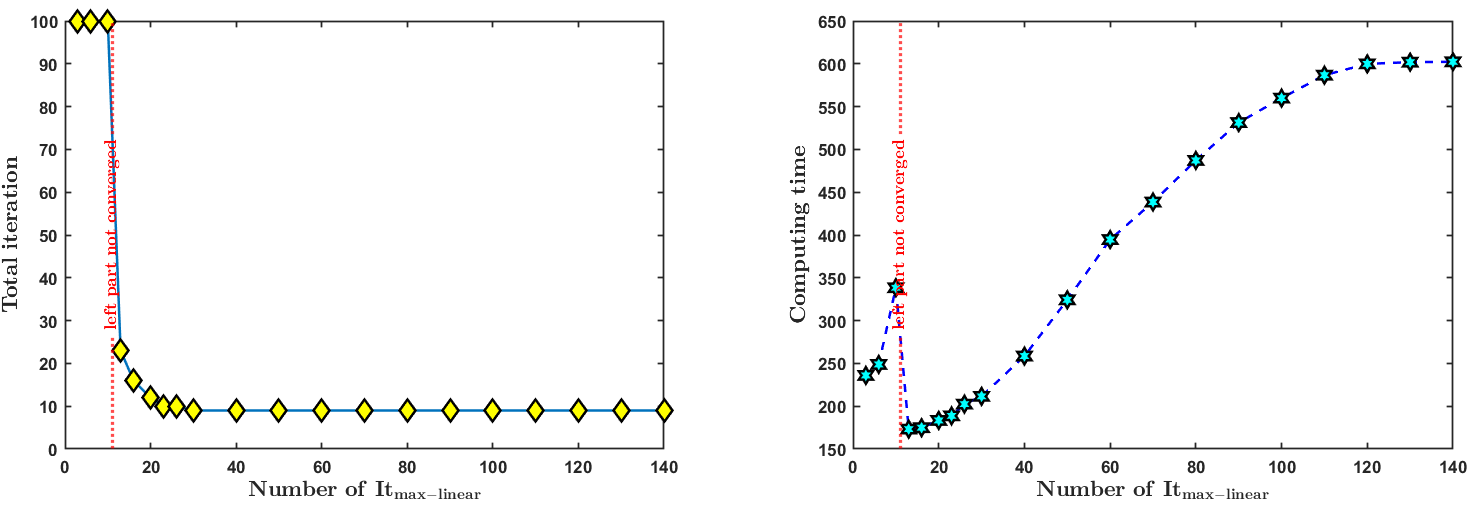}
	\caption{The outer iterations (left) and CPU time (right) of the enhanced ParaSLRF with different $\rm It_{\text{max-linear}}$ for the pyramid S1 model ($N=36339$).}
	\label{fig: inner_pyramid}
\end{figure}

\begin{figure}[htbp]
	\centering
	\includegraphics[scale=0.4]{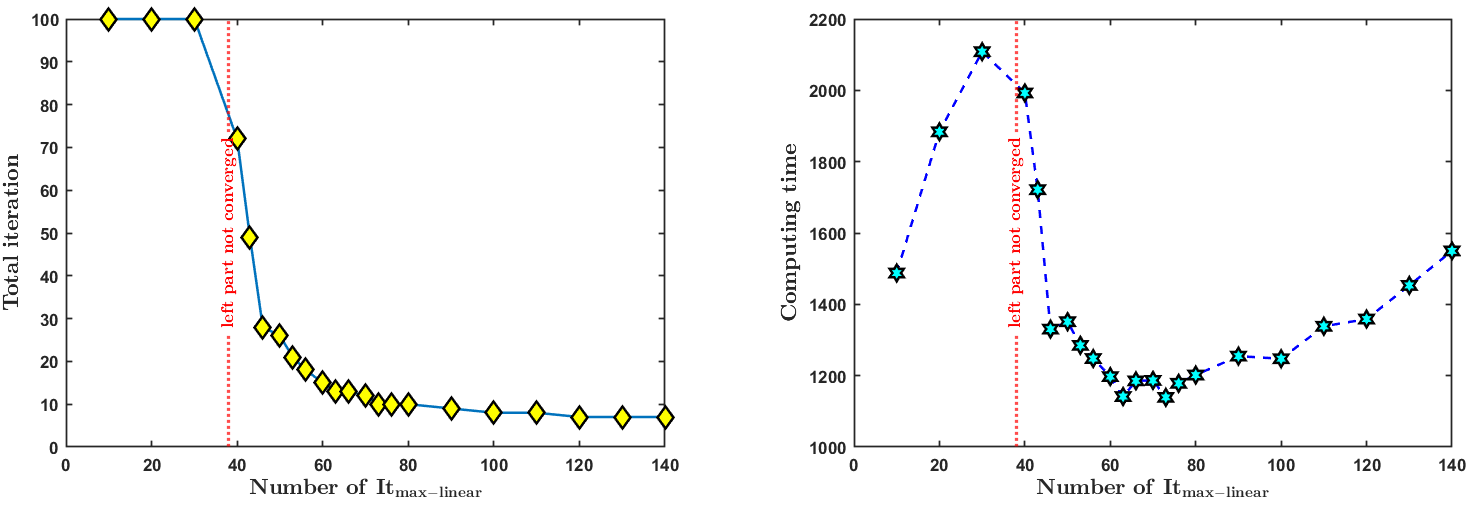}
	\caption{The outer iterations (left) and CPU time (right) of the enhanced ParaSLRF with different $\rm It_{\text{max-linear}}$ for hollow platform model ($N=37161$).}
	\label{fig: inner_trunc}
\end{figure}

\begin{figure}[htbp]
	\centering
	\includegraphics[scale=0.4]{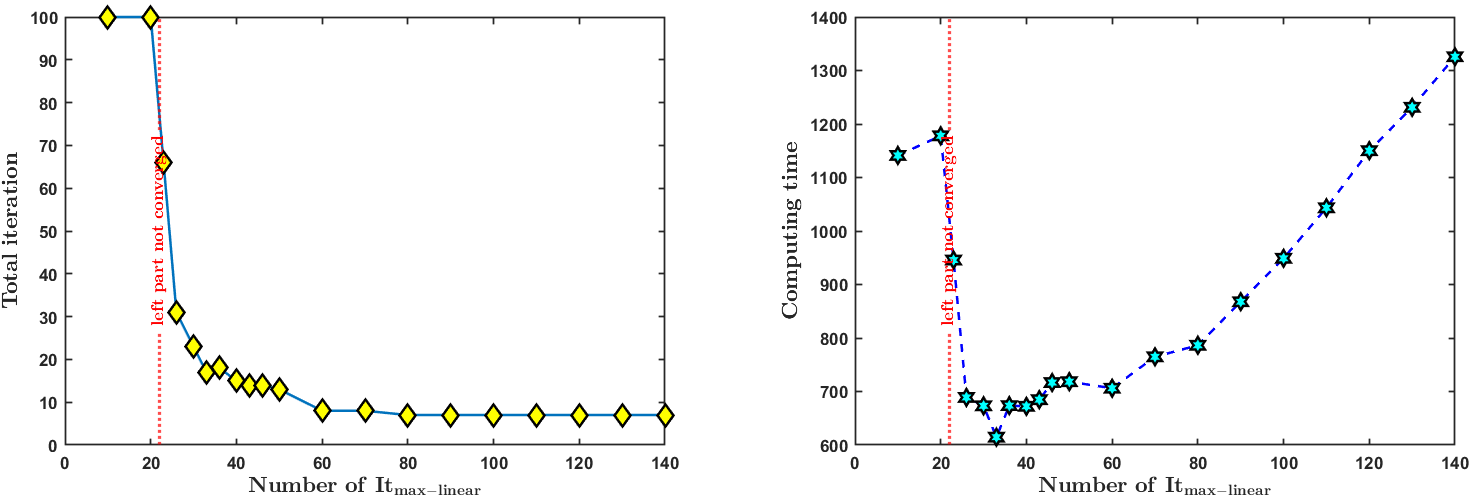}
	\caption{The outer iterations (left) and CPU time (right) of the enhanced ParaSLRF with different $\rm It_{\text{max-linear}}$ for fish-like model ($N=55248$).}
	\label{fig: inner_hollow}
\end{figure}

Although the new algorithm adds overhead due to orthogonalization, it significantly reduces the cost of the linear solves, compared to Algorithm~\ref{alg:Parallel Rational filter method}.
It is not easy to set a default value for the maximum number of inner linear iterations $\rm It_{\text{max-linear}}$.
In order to investigate the effect of $\rm It_{\text{max-linear}}$ on the convergence and performance of the enhanced ParaSLRF,
we computed the first 100 eigenpairs of aforementioned models, for different values of $\rm It_{\text{max-linear}}$. 
The CPU time and the number of outer iterations for each task can be found in the Figures~\ref{fig: inner_pyramid}---\ref{fig: inner_hollow}, respectively.

In Figures~\ref{fig: inner_pyramid}---\ref{fig: inner_hollow} the red dashed line shows the minimal value of $\rm It_{\text{max-linear}}$ to obtain convergence of the wanted eigenpairs in 100 outer iterations.
It is not surprising that the best value of $\rm It_{\text{max-linear}}$ for achieving the lowest CPU time is problem dependent.
Although a high $\rm It_{\text{max-linear}}$ can reduce the number of outer iterations this may also lead to excessive CPU time. Therefore, a careful trade-off is required when selecting this parameter.
Through experimental validation, the optimal 
$\rm It_{\text{max-linear}}$ was set to 13, 63, and 33 for the pyramid S1, hollow platform, and fish-like models, respectively.

To demonstrate the improved performance, 
we compute the first 20 and 100 smallest eigenpairs of aforementioned mechanical models with corresponding best value of $\rm It_{\text{max-linear}}$.
The computational environment and the notation used for reporting the results are the same as previously described in Section~\ref{subsect: results}.
The results are shown in  Table~\ref{tab: enhanced results}. 
\begin{table}[htbp]
  \centering
  \caption{Performance improvement.}
  \begin{threeparttable}
    \begin{tabular}{c|cccccc}
    \multicolumn{1}{c}{Model} & Interval & Method & Iter  & Time$_{\text{Total}}$  & Time$_{\text{MPIWait}}^{\text{Max}}$ & Prop (\%) \\
    \midrule
    \multicolumn{1}{c|}{\multirow{4}[4]{*}{Pyramid S1}} & \multicolumn{1}{c}{\multirow{2}[2]{*}{\shortstack{NEV=20 \\ (0, 1531]}}} & ParaSLRF &  10   & 320.852      &   8.175   &  2.50 \\
          &       & Enhanced ParaSLRF &  20    &   35.247  & 0.514 & 1.46\\
\cmidrule{2-7}          & \multicolumn{1}{c}{\multirow{2}[2]{*}{\shortstack{NEV=100 \\(0, 4273]}}} & ParaSLRF &   9   &   1153.512    &     10.824  &  0.93 \\
          &       & Enhanced ParaSLRF &   23    &   173.227 &  1.347
  & 0.78 \\
    \midrule
    \multirow{4}[4]{*}{\shortstack{Hollow \\ platform}} & \multicolumn{1}{c}{\multirow{2}[2]{*}{\shortstack{NEV=20 \\ (0, 215]}}} & ParaSLRF &   6   &  1540.356     &    34.959   &  2.26 \\
          &       & Enhanced ParaSLRF &   15    &   244.410 &  3.359 & 1.37\\
\cmidrule{2-7}          & \multicolumn{1}{c}{\multirow{2}[2]{*}{\shortstack{NEV=100 \\ (0, 1507]}}} & ParaSLRF &  7   &   5426.306    &      124.050 & 2.28  \\
          &       & Enhanced ParaSLRF &    13   &   1147.207 & 16.263 & 1.42 \\
    \midrule
    \multirow{4}[4]{*}{Fish-like} & \multicolumn{1}{c}{\multirow{2}[2]{*}{\shortstack{NEV=20 \\ (0, 1858]}}} & ParaSLRF &  7    &   1607.061    &   88.895    &  5.53 \\
          &       & Enhanced ParaSLRF &  50 &  258.127 & 4.319&  1.67   \\
\cmidrule{2-7}          & \multicolumn{1}{c}{\multirow{2}[2]{*}{\shortstack{NEV=100 \\ (0, 10239]}}} & ParaSLRF &    7   &      3838.613 &   147.550    &  3.84 \\
          &       & Enhanced ParaSLRF &  17     &   613.464 & 8.967 & 1.46 \\
    \bottomrule
    \end{tabular}%
    \end{threeparttable}
  \label{tab: enhanced results}%
\end{table}%

From Table~\ref{tab: enhanced results}, we can observe that the computational  cost and CPU time consumption are significantly reduced although the two applied strategies slightly increase the number of outer iterations. 
In addition, the enhanced version is at least 4.7 times faster than the original one for all tasks and can be up to 9.1 times faster when computing 20 eigenpairs of pyramid S1 model.

We also find that, for the enhanced ParaSLRF, the extra MPI waiting time reduces significantly; the proportion of the extra MPI waiting time to the total CPU time also reduces further, which represents better utilization of resources and load-balancing.

\subsection{Strong scalability}

To investigate the strong scalability of the enhanced ParaSLRF, we use the pyramid S2 model represented in~\S\ref{subsect: pyramid model} with refined mesh, as test example.
The discretized system has 268,353  degrees of freedom. 
We set $\rm It_{\text{max-linear}}=60$ and intend to compute the first 100 smallest eigenpairs of this system, with different computing resources configuration. 
The computations were carried out on a cluster where each computational node has 2 Intel Xeon Gold 6240 CPUs@2.6 GHz (18 cores per CPU), and 192 GiB RAM.

\begin{figure}[htbp]
\flushleft
\includegraphics[scale=0.55]{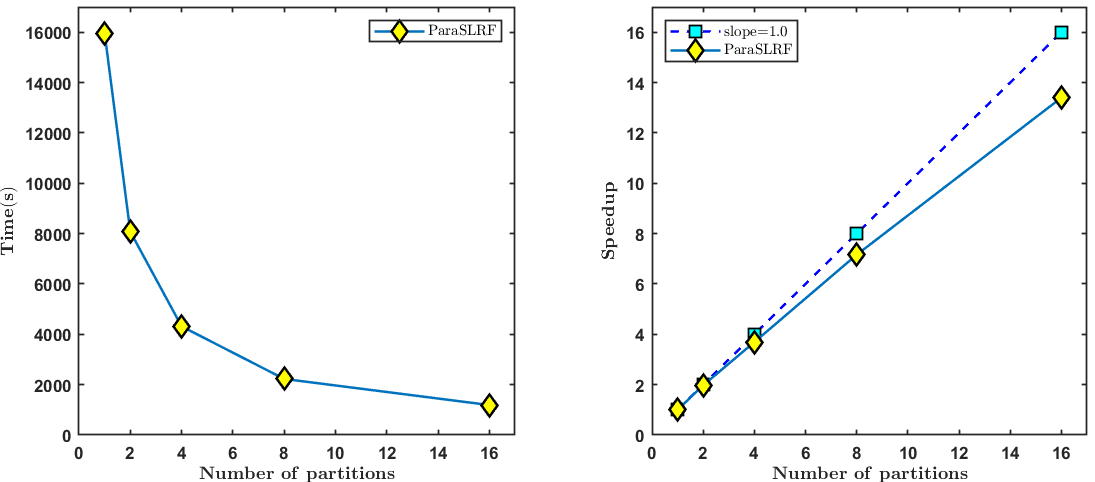}
\caption{Strong scalability of the Level-1 parallelism.}
\label{fig: L1_scalability}
\end{figure}

For the Level-1 scalability test, we use $N=16$ poles and split the poles into $n_{part}=1,2,4,8,16$ parts, respectively. We allocate $np_{sub}=18$ processes bound to one socket for each part. 
The total number of processes $np$ varies as 18,36,72,144,288.
We record the CPU time and compute the corresponding speed-up. The results are shown in Figure~\ref{fig: L1_scalability}. 

From Figure~\ref{fig: L1_scalability} we observe that the enhanced ParaSLRF performs well on the Level-1 strong scalability, which is close to  linear.
These results indicate that prioritizing resources allocation at the Level-1 is critical for achieving a good parallel performance.

\begin{figure}[htbp]
\flushleft
\includegraphics[scale=0.55]{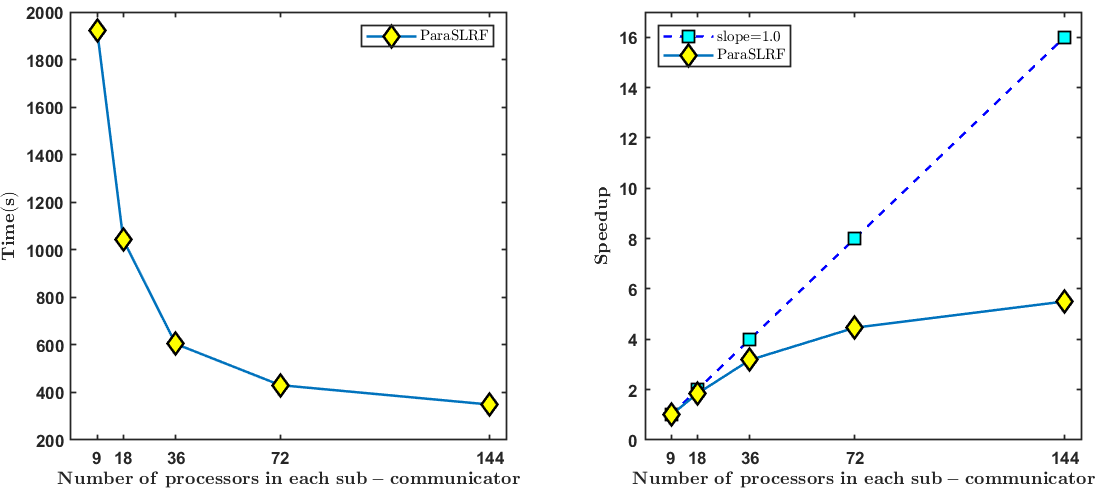}
\caption{Strong scalability of the Level-2 parallelism.}
\label{fig: L2_scalability}
\end{figure}

As for the Level-2 scalability test, it is actually the strong scalability test of the linear iterative solver.
Thus, we do not need to use a large number of poles for this test. Here, we only use $N=4$ poles and split all processes into $n_{part}=4$ parts.
The number of allocated processes $np_{sub}$ starts from $9,18,36,72,144$ for each sub-communicator (with a total of 576 cores).
The CPU time and corresponding speed-up for each computation task are represented in Figure~\ref{fig: L2_scalability}.

From Figure~\ref{fig: L2_scalability} we observe that the Level-2 strong scalability of the enhanced ParaSLRF is not so good. On the one hand, solving sparse linear systems is memory-bound, which means that it is hard to achieve the ideal scalability for any linear solver; and on the other hand, as the number of processes grows rapidly the number of degrees of freedom handled by each process is smaller, which can pull down the performance of GAMG.
In general, this is a normal scalability performance of the GCR iterative method preconditioned by GAMG.

\subsection{Fixed computing resources}
\label{subsect: fixed resources}
Because of the two levels of parallelism within the enhanced ParaSLRF, it is valuable to explore strategies for making the most efficient use of fixed computing resources (e.g., the limited number of processes), to achieve the best parallel performance.
We have to consider how to allocate computing resources in each parallel level.
If 4 computational nodes (with a total of 144 cores) are available and 16 poles are used for the computations, there are five ways to allocate the resources across the two-level parallelism.
It is reasonable to set $n_{part}=1,2,4,8,16$ for the number of sub-communicators in Level-1 and utilize $np_{sub}=144,72,36,18,9$ processes in each sub-communicator, respectively.

Table~\ref{tab: fixed resources} shows the CPU time for solving the first 100 eigenpairs of the pyramid S2 model, with different distributions of computing resources across the Level-1 and Level-2. 
It is no surprise that it is optimal to place as many resources as possible at the Level-1 and the best performance is achieved when the required linear systems only involve one pole in each sub-communicator.

\begin{table}[htbp]
  \centering
  \caption{The CPU time for computing 100 eigenpairs of the Pyramid S2 model with different distributions of computing resources. Increasing the number of MPI processes at Level-2 (i.e., $np_{sub}$) reduces the memory required by each process (i.e., the number of rows of the vectors and matrices assigned to each process).}
    \begin{tabular}{cccc}
    \toprule
    $n_{part}$ & $np_{sub}$ & rows  & Time(s) \\
    \midrule
    1     & 144    & 1864  & 3420.508 \\
    2     & 72    & 3727 & 2817.676 \\
    4     & 36    & 7454 & 2369.271 \\
    8     & 18     & 14909 & 2199.769 \\
    16    & 9      & 29817  &  2128.285 \\
    \bottomrule
    \end{tabular}%
  \label{tab: fixed resources}%
\end{table}%

\section{Conclusion}
\label{sect:conclusion}
In this paper, we developed a parallel shifted Laplace rational filter method (ParaSLRF) for computing all eigenvalues located in the interval $(0,\gamma]$, and corresponding eigenvectors, for generalized eigenvalue problems arising from mechanical vibrations.

Numerical results indicate that, compared to the rational filters based on quadrature rules, ParaSLRF performs much better in terms of computational cost and CPU time. 
Moreover, we observe that there is less communication waiting overhead between sub-communicators for ParaSLRF, i.e.,  it offers excellent load-balancing capability which makes ParaSLRF highly suitable for parallel computing.

In addition, we introduced two improvements to enhance the performance of ParaSLRF. First, converged eigenpairs are locked by making the filtered vectors $\matB$-orthogonal to the converged vectors. A significant reduction of the computational cost is obtained with only a small additional overhead.
Second, we selected a good initial guess for the linear iterative solver and fixed the maximum linear iteration number to reduce the computational cost further. This led to a better load-balancing.
Scalability tests show that the enhanced version has very good two-level parallel scalability.
Moreover, we observe that distributing as much computing resources as possible at Level-1 is the optimal distribution with limited resources.


 \bibliographystyle{elsarticle-num} 
 \bibliography{rationalfilter}





\end{document}